\documentstyle [12pt,epsf, graphicx, amssymb]{article}

\begin{document}
\def\l{\lambda}
\def\m{\mu}
\def\a{\alpha}
\def\b{\beta}
\def\g{\gamma}
\def\d{\delta}
\def\e{\epsilon}
\def\o{\omega}
\def\O{\Omega}
\def\v{\varphi}
\def\t{\theta}
\def\r{\rho}
\def\bs{$\blacksquare$}
\def\bp{\begin{proposition}}
\def\ep{\end{proposition}}
\def\bt{\begin{th}}
\def\et{\end{th}}
\def\be{\begin{equation}}
\def\ee{\end{equation}}
\def\bl{\begin{lemma}}
\def\el{\end{lemma}}
\def\bc{\begin{corollary}}
\def\ec{\end{corollary}}
\def\pr{\noindent{\bf Proof: }}
\def\note{\noindent{\bf Note. }}
\def\bd{\begin{definition}}
\def\ed{\end{definition}}
\def\C{{\mathbb C}}
\def\P{{\mathbb P}}
\def\Z{{\mathbb Z}}
\def\d{{\rm d}}
\def\deg{{\rm deg\,}}
\def\deg{{\rm deg\,}}
\def\arg{{\rm arg\,}}
\def\min{{\rm min\,}}
\def\max{{\rm max\,}}

\newcommand{\norm}[1]{\left\Vert#1\right\Vert}
\newcommand{\abs}[1]{\left\vert#1\right\vert}

\newcommand{\set}[1]{\left\{#1\right\}}
\newcommand{\setb}[2]{ \left\{#1 \ \Big| \ #2 \right\} }

\newcommand{\IP}[1]{\left<#1\right>}
\newcommand{\Bracket}[1]{\left[#1\right]}
\newcommand{\Soger}[1]{\left(#1\right)}

\newcommand{\Integer}{\mathbb{Z}}
\newcommand{\Rational}{\mathbb{Q}}
\newcommand{\Real}{\mathbb{R}}
\newcommand{\Complex}{\mathbb{C}}

\newcommand{\eps}{\varepsilon}
\newcommand{\To}{\longrightarrow}
\newcommand{\varchi}{\raisebox{2pt}{$\chi$}}

\newcommand{\E}{\mathbf{E}}
\newcommand{\Var}{\mathrm{var}}

\def\squareforqed{\hbox{\rlap{$\sqcap$}$\sqcup$}}
\def\qed{\ifmmode\squareforqed\else{\unskip\nobreak\hfil
\penalty50\hskip1em\null\nobreak\hfil\squareforqed
\parfillskip=0pt\finalhyphendemerits=0\endgraf}\fi}

\renewcommand{\th}{^{\mathrm{th}}}
\newcommand{\Dif}{\mathrm{D_{if}}}
\newcommand{\Difp}{\mathrm{D^p_{if}}}
\newcommand{\GHF}{\mathrm{G_{HF}}}
\newcommand{\GHFP}{\mathrm{G^p_{HF}}}
\newcommand{\f}{\mathrm{f}}
\newcommand{\fgh}{\mathrm{f_{gh}}}
\newcommand{\T}{\mathrm{T}}
\newcommand{\K}{^\mathrm{K}}
\newcommand{\PghK}{\mathrm{P^K_{f_{gh}}}}
\newcommand{\Dig}{\mathrm{D_{ig}}}
\newcommand{\for}{\mathrm{for}}
\newcommand{\End}{\mathrm{end}}

\newtheorem{th}{Theorem}[section]
\newtheorem{lemma}{Lemma}[section]
\newtheorem{definition}{Definition}[section]
\newtheorem{corollary}{Corollary}
\newtheorem{proposition}{Proposition}[section]

\begin{titlepage}

\begin{center}

\topskip 5mm

{\LARGE{\bf {Remez-Type Inequality for Discrete Sets}}} \vskip 8mm

{\large {\bf Y. Yomdin}}

\vspace{6 mm}
\end{center}

{Department of Mathematics, The Weizmann Institute of Science,
Rehovot 76100, Israel}

\vspace{6 mm}

{e-mail: yosef.yomdin@weizmann.ac.il}

\vspace{1 mm}

\vspace{6 mm}
\begin{center}

{ \bf Abstract}
\end{center}

{\small The classical Remez inequality bounds the maximum of the
absolute value of a polynomial $P(x)$ of degree $d$ on $[-1,1]$
through the maximum of its absolute value on any subset $Z$ of
positive measure in $[-1,1]$. Similarly, in several variables the
maximum of the absolute value of a polynomial $P(x)$ of degree $d$
on the unit cube $Q^n_1 \subset {\mathbb R}^n$ can be bounded
through the maximum of its absolute value on any subset $Z\subset
Q^n_1$ of positive $n$-measure.

The main result of this paper is that the $n$-measure in the Remez
inequality can be replaced by a certain geometric invariant
$\omega_d(Z)$ which can be effectively estimated in terms of the
metric entropy of $Z$ and which may be nonzero for discrete and
even finite sets $Z$.}

\end{titlepage}

\newpage

\section{Introduction}
\setcounter{equation}{0}

The classical Remez inequality (\cite{Rem}) bounds the maximum of
the absolute value of a polynomial $P(x)$ of degree $d$ on
$[-1,1]$ through the maximum of its absolute value on any subset
$Z$ of positive measure in $[-1,1]$. More accurately:

\medskip

\noindent {\it Let $P(x)$ be a polynomial of degree $d$. Then for
any measurable $Z\subset [-1,1]$ \be \max_{[-1,1]} \vert P(x)
\vert \leq T_d({{4-\mu}\over {\mu}})\max_Z \vert P(x) \vert ,\ee
where $\mu=\mu_1(Z)$ is the Lebesgue measure of $Z$ and
$T_d(x)=cos(d \ arccos(x))$ is the $d$-th Chebyshev polynomial}.

\medskip

In several variables the generalization of (1.1) was obtained in
\cite{Bru.Gan}:

\bt Let ${\cal
B}\subset {\mathbb R}^n$ be a convex body and let $\O\subset {\cal
B}$ be a measurable set. Then for any real polynomial
$P(x)=P(x_1,\dots,x_n)$ of degree $d$ we have \be \sup_{{\cal
B}}\vert P \vert \leq T_d ({{1+(1-\lambda)^{1\over n}}\over
{1-(1-\lambda)^{1\over n}}}) \sup_{\O}\vert P \vert.\ee Here
$\lambda= {{\mu_n(\O)}\over {\mu_n({\cal B})}},$ with $\mu_n$ being
the Lebesgue measure on ${\mathbb R}^n$. This inequality
is sharp and for $n=1$ it coincides with the classical Remez
inequality.\et

It is well known that the inequality of the form (1.1) or (1.2)
may be true also for some sets $Z$ of measure zero and even for
certain discrete or finite sets $Z$. Let us mention here only a
couple of the most relevant results in this direction: in
\cite{Ber2,Ber3,Cop.Riv,Nad,Rah,Yom.Zah,Zah} such inequalities are
provided for $Z$ being a regular grid in $[-1,1]$ (in \cite{Ber3}
trigonometric polynomials are considered). Let us mention
\cite{Ber1} where a ``dual" problem is considered of interpolation
by polynomials of degree higher than the number of the nodes. In
\cite{Fav} discrete sets $Z\subset [-1,1]$ are studied (see
Section 2.2 below). An invariant $\phi_Z(d)$ is defined and
estimated in some examples, which is the best constant in the
Remez-type inequality of degree $d$ for the couple $(Z\subset
[-1,1])$. Below in Definition 1.1 we call this invariant (extended
to any dimension) the Remez $d$-span of $Z$ and denote it
$R_d(Z)$.

On the other hand, recently in \cite{Bru.Bru2} Remez inequality
has been extended (for {\it complex} polynomials of $n$ variables)
to subsets $Z$ of positive Hausdorff $s$-measure, $s> 2n-2$. Here
$2n-2$ is the real dimension of a zero set of such a polynomial,
so the result has a natural geometric interpretation: Remez-type
inequalities are true for $Z$ having Hausdorff dimension larger
than the dimension of the corresponding zero sets. For real
polynomials of $n$ variables, under the above assumption on $Z$,
an integral version of the Remez inequality was proved in
\cite{Bru.Bru2}, and a question was posed of the existence of a ``strong"
Remez-type inequality (of the form (1.2)).

In \cite{Bru} estimates have been obtained for covering numbers of sub-level sets of
families of analytic functions depending analytically on a parameter. Using these
estimates strong Remez type inequalities have been proved for the restrictions of analytic
functions to certain fractal sets. The existence of such inequalities was conjectured
in \cite{Bru.Bru2}.

In \cite{Bor.Naz.Sod}, \cite{Naz.Sod.Vol1}-\cite{Naz.Sod.Vol3} analytic and
quasi-analytic functions have been studied from a similar point of view.

\smallskip

For one complex variable results similar to Remez inequality are
provided by the classical Cartan lemma
(see, for example, \cite{Gor,Erd} and references therein):

\smallskip

\noindent {\it Let $P(z)$ be a monic polynomial of one complex
variable of degree $d$. For any given $\e>0$ consider
$V_{\e^d}(P)=\{z\in {\mathbb C}, \ \vert P(z) \vert \leq \e^d\}.$
Then $V_{\e^d}(P)$ can be covered by at most $d$ complex discs $D_j$ with
radii $r_j, \ j=1,\dots,d$ such that $\sum_{j=1}^d r_j \leq 2e
\e.$}

\smallskip

In \cite{Zer1} (see also \cite{Zer2}) a generalization of the
Cartan lemma to plurisubharmonic functions was obtained
which leads, in particular, to the bounds on the size of sub-level
sets similar to those obtained in \cite{Bru.Bru2}.

\medskip

In the present paper we would like to address a general problem of
characterizing sets $Z$ for which Remez-type inequality is valid.
Having this in mind let us give the following definition: \bd A
set $Z\subset Q^n_1 \subset {\mathbb R}^n$ is called $d$-definite
if any real polynomial $P(x)=P(x_1,\dots,x_n)$ of degree $d$ bounded in
absolute value by $1$ on $Z$ is bounded in absolute value by a
certain constant $C(d,Z)$ (not depending on $P$) on $Q^n_1$. The
minimum $R_d(Z)$ of all such constants $C(d,Z)$ is called the
Remez $d$-span of $Z$.\ed In view of the above-mentioned results
the following problem looks natural and important:

\medskip

{\it Characterize (through their metric geometry) all the sets
$Z\subset Q^n_1 \subset {\mathbb R}^n$ with the finite Remez
$d$-span $R_d(Z)$ and compute $R_d(Z)$ for such $Z$ in ``geometric" terms.}

\medskip

In principle,
there is a very simple answer to this question: {\it $R_d(Z)=\infty$ if and
only if $Z$ is contained in a zero set of a certain nonzero
polynomial $P$ of degree $d$}. Indeed, in the opposite case $\sup_{Z}\vert P \vert$
and $\sup_{Q_1^n}\vert P \vert$ both are norms on the finite dimensional space
of polynomials of degree $d$. However, in general it is not easy to reformulate this
condition in ``effective geometric terms" and to provide explicit
bounds on $R_d(Z)$ starting with an explicitly given $Z$. For finite
sets $Z$ it is possible (in principle) to write an explicit answer
through the ``interpolation systems" (see, for example, \cite{Hil,Lor,For}
and references therein). But to analyze, for instance, the asymptotic
behavior of $R_d(Z)$ as $d\rightarrow \infty$ for a ``fractal" $Z$ may
be a tough problem (compare \cite{Fav} and Section 2.2 below).

\medskip

In the present paper we construct for subsets $Z\subset Q^n_1$ a
simple geometric invariant $\omega_d(Z)$ which we call a metric
$d$-span of $Z$. The metric $d$-span $\o_d(Z)$ can be effectively
estimated in terms of the metric entropy of $Z$ and it may be
nonzero for discrete and even finite sets $Z$. Our main result is
that {\it the $n$-measure in the Remez-type inequality (1.2) can be replaced
by $\omega_d(Z)$}.

\medskip

To define $\omega_d(Z)$ let us recall that the covering number $M(\e,A)$
of a metric space $A$ is the minimal number of closed $\e$-balls covering $A$ (see
\cite{Fal,Kol,Kol.Tih,Min}).
Below $A$ will be subsets of ${\mathbb R}^n$ equipped with the $l^\infty$ metric.
So the $\e$-balls in this metric are the cubes $Q^n_\e$.

For a polynomial $P$ on ${\mathbb R}^n$ let us consider the
sub-level set $V_\rho(P)$ defined by $V_\rho(P)=\{x\in Q^n_1,
\vert P(x) \vert \leq \rho \}$. The following result is provided
by (\cite{Vit1,Vit2,Iva}): \bt ({\bf Vitushkin's bound}) For
$V=V_\rho(P)$ as above \be M(\e,V)\leq
\sum_{i=0}^{n-1}C_i(n,d)({1\over \e})^i+
\mu_n(V)({1\over{\e}})^n,\ee with $C_i(n,d)=C'_i(n)(2d)^{(n-i)}$.
For $n=1$ we have $M(\e,V)\leq d + \mu_1(V)({1\over{\e}}),$ and
for $n=2$ we have
$$M(\e,V)\leq (2d-1)^2 + 8d({1\over \e}) + \mu_2(V)({1\over{\e}})^2.$$\et
For $\e>0$ we denote by $M_{n,d}(\e)$ (or shortly $M_d(\e)$) the
polynomial of degree $n-1$ in ${1\over \e}$ as appears in (1.3):
\be M_d(\e)=\sum_{i=0}^{n-1}C_i(n,d)({1\over \e})^i.\ee In
particular,
$$M_{1,d}(\e)=d, \ M_{2,d}(\e)=(2d-1)^2 + 8d({1\over \e}).$$ Now for each
subset $Z\subset Q_1^n$ (possibly discrete or finite) we introduce
the metric $(n,d)$-span of $Z$ via the following definition: \bd
Let $Z$ be a subset in $Q^n_1\subset {\mathbb R}^n$. Then the
metric $(n,d)$-span (or shortly $d$-span) $\omega_d(Z)$ is defined
as \be \omega_d(Z) = \sup_{\e>0} {\e}^n[M(\e,Z)- M_d(\e)]. \ee \ed
Now we are ready to state the main result of this paper which, in
particular, provides a partial answer to the above-stated general
problem of describing sets with finite Remez $d$-span: \bt If
$\omega_d(Z) = \o >0$ then $R_d(Z)$ is finite and satisfies \be
R_d(Z)\leq R(\o):= T_d ({{1+(1-\o)^{1\over n}}\over
{1-(1-\o)^{1\over n}}}),\ee where $T_d(x)=cos(d \arccos(x))$ is
the $d$-th Chebyshev polynomial.\et This theorem is proved in
Section 2 via a combination of inequalities (1.2) and (1.3), the
last being reinterpreted as in Theorem 2.1 below. The bound in
(1.6) is finite if and only if $\o > 0$. As an immediate corollary
of Definition 1.2 and Theorem 1.3 we obtain the following general
condition for positivity of $\o_d(Z)$: \bc The $d$-span $\o_d(Z)$
is positive if and only if for certain $\e>0$ we have $M(\e,Z)>
M_d(\e)$. \ec

\smallskip

Corollary 1 establishes finiteness of $R_d(Z)$ for a large class
of sets. First of all, this is true for $Z$ having a Hausdorff
dimension $\dim_H(Z)$ greater than $n-1$. Together with an
explicit bound given in Section 3 below this provides a partial
answer to the question posed in \cite{Bru.Bru2} of existence in
this case of a Remez-type inequality of the form (1.2).

In fact, we can replace the Hausdorff dimension by the entropy (or
the box) dimension $\dim_e(Z)$. The entropy dimension is defined
in terms of the asymptotic behavior of the covering number
$M(\e,Z)$ as $\e$ tends to $0$. It is larger (and often strictly
larger) than the Hausdorff dimension and it may take any value up
to $n$ also for countable subsets of ${\mathbb R}^n$.

\smallskip

However, Corollary 1 shows finiteness of $R_d(Z)$ also for
sufficiently large sets of dimension exactly $n-1$: \bc Let $Z$ be
a $C^1$-hypersurface in $Q^n_1$ with $\mu_{n-1}(Z)> C_{n-1}(n,d)$.
Then $Z$ is $d$-definite.\ec By the virtue of the covering number
{\it also sufficiently dense finite subsets of $d$-definite sets
are themselves $d$-definite.} Thus Corollary 1 provides a class of
examples of finite $d$-definite subsets $Z\subset B^n_1$: roughly,
those are sufficiently dense finite subsets of sets of dimension
$n-1$ or higher.

It is important to analyze the behavior of the Remez $d$-span
$R_d(Z)$ for finite (and general) sets $Z$ in terms of their
metric structure. Here an appropriate invariant
may be the so-called $\beta$-spread, introduced in
\cite{Yom1,Yom2}. A very closely related notion is the
``$\beta$-weight of a minimal spanning trees" (see
\cite{Koz.Lot.Stu} and references therein). Some initial results
in this direction are given in Section 3.2.6 below.

\medskip

Let us stress that {\it the sufficient condition for a set $Z$ to
be $d$-definite provided by Corollary 1 is not necessary in
general}: any small piece $W'$ of an irreducible algebraic
hypersurface $W$ in ${\mathbb R}^n$ of degree $d_1>d$ is
$d$-definite, since it is not contained in any algebraic
hypersurface of degree $d$. The same is true for transcendental
hypersurfaces $W$ as well as for transcendental (or algebraic of
high degree) sets $W$ of smaller dimensions, in particular, for
curves.

On the other hand, in each of these situations, if we take a
sufficiently small piece $W'$ of $W$ inside the unit cube, the
$(n-1)$-area $\mu_{n-1}(W')$ of $W'$ is much less than
$C_{n-1}(n,d)$. Consequently, $M(\e,W')\asymp
\mu_{n-1}(W')({1\over \e})^{n-1}$ is always strictly less than
$C_{n-1}(n,d)({1\over \e})^{n-1}$. Therefore by Corollary 1 we
have $\o_d(W')=0$.

\medskip

This stresses once more the importance (and apparent difficulty)
of the problem of a ``geometric characterization" of the
$d$-definite sets.

\smallskip

The behavior of polynomials on discrete sets
plays an important role in the Whitney problem of extension of
differentiable functions from closed sets (\cite{Whi}). In
particular, there is an apparent relation of the Remez type
inequalities with the problem of extending ``finite differences"
to higher dimensions.  See \cite{Bie.Mil.Paw,Bru.Shv,Fef} for some
results representing recent progress in the Whitney problem.

\medskip
\medskip

The paper is organized as follows: in Section 2 we prove Theorem
1.3 and state (in a simplified form) some of the main results
describing the behavior of the $d$-span in specific situations.
The proofs of these results are postponed till Section 3, where
they are given with all the required details and accurate (but
somewhat cumbersome) constants.

\medskip

The author would like to thank A. Brudnyi, V. Katsnelson, B.
Nadler and M. Sodin for useful discussions and for providing
important references, and the referee for suggesting significant
improvements of the presentation.

\section{Proof of Theorem 1.3 and some basic examples of $d$-definite sets}
\setcounter{equation}{0}

The following theorem relates the $n$-volume of the sub-level sets
$V_\rho$ of a polynomial of degree $d$ with the metric
$(d,n)$-span $\o_d(Z)$ of subsets $Z\subset V_\rho$: \bt Let
$P(x)=P(x_1,\dots,x_n)$ be a polynomial of degree $d$ and let
$Z\subset Q^n_1$ be a given set. Then if $Z\subset V_\rho(P)$ for
a certain $\rho\geq 0$ then we have
$$\mu_n(V_\rho(P))\geq \omega_d(Z),$$ where $\mu_n$ denotes, as
above, the Lebesgue $n$-measure. \et \pr This fact follows
directly from the Vitushkin bound on the covering number of the sub-level
sets $V_\rho$ given in Theorem 1.2 above: for any polynomial
$P(x)=P(x_1,\dots,x_n)$ of degree $d$, for $Z\subset V_\rho(P)$, and for
any $\e>0$ we have
$$ M(\e,Z)\leq M(\e,V_\rho(P)) \leq M_d(\e) + \mu_n(V_\rho(P)) ({1\over{\e}})^n.$$
Consequently, for any $\e>0$ we have \be \mu_n(V_\rho(P)) \geq
\e^n[M(\e,Z)- M_d(\e)],\ee and we can take a supremum with respect
to $\e$. This completes the proof. $\blacksquare$

\medskip

\noindent{\bf Proof of Theorem 1.3} Assume that $P$ is bounded in absolute value
by $1$ on $Z$. Then we have $Z\subset V_1(P)$. By Theorem 2.1
$\mu_n(V_1(P))\geq \o_d(Z)=\o$. Now since $P$ is bounded in
absolute value by $1$ on $V_1(P)$ by definition, we can apply the
Yu. Brudnyi-Ganzburg inequality (Theorem 1.1 above)
with ${\cal B}=Q^n_1$ and $\O=V_1(P)$. This completes the proof.
$\blacksquare$

\smallskip

Let us now study some specific classes of $Z$. \bt A set $Z\subset
Q^n_1$ of positive $s$-Hausdorff measure, $s>n-1$, is $d$-definite
for any $d$.\et \pr This follows directly from Theorem 3.2,
Section 3 below, where also a lower bound for $\o_d(Z)$ is given.
$\blacksquare$

\medskip

The invariant $\o_d$ is strong enough to prove that sets of
dimension exactly $n-1$ are definite, assuming their
$(n-1)$-Hausdorff measure is big enough. \bt A set $Z\subset
Q^n_1$ with $(n-1)$-Hausdorff measure $H_{n-1}(Z)$ satisfying
$$H_{n-1}(Z)>2 \sqrt{n^{n-1}}C_{n-1}(n,d)$$ is $d$-definite. In
particular, any curve $Z\subset Q^2_1$ of the length $l(Z)$
satisfying $l(Z)>16\sqrt 2 d$ is $d$-definite. \et \pr This
follows directly from Corollary 4, Section 3 below, where also a
lower bound for $\o_d(Z)$ is given. $\blacksquare$

\subsection{Bounding Remez $d$-span via Minimal spanning trees}

Let us now consider finite sets $Z$. By virtue of the definitions
any sufficiently dense finite subset of a set with positive
$d$-span also has positive $d$-span. See Theorem 3.4 of Section 3
below specifying the choice of such a dense finite subset in each
of the cases considered above. One result addressing the specific
geometry of $Z$ is the following (the distance below is with
respect to the $l^\infty$-norm on ${\mathbb R}^n$):

\bt Let a degree $d$ and a finite subset $Z\subset Q^n_1$ be
given, and let $\e_0$ be the minimal distance between the points
of $Z$. Assume that $\vert Z \vert = p > M_d(\e_0).$ Then the set
$Z$ is $d$-definite. In particular, any set with more than $d$
points in $Q^1_1$ is $d$-definite. Any set $Z\subset Q^2_1$ with
the number of points larger than $(2d-1)^2+8d({1\over \e_0})$ is
$d$-definite. \et \pr This follows from Definition 1.2
since $M(\e_o,Z)=\vert Z \vert = p.$ $\blacksquare$

\medskip

\noindent{\bf Remark} Let us stress the importance of the assumption
$Z\subset Q^n_1$ in Theorem 2.4. Without it we could take all the points
of $Z$ on the same straight line. However, inside the cube $Q^n_1$ the
points of $Z$ must form an ``essentially $n$-dimensional" configuration
in order to satisfy the inequality $\vert Z \vert > M_d(\e_0)$.

\medskip

Following the direction of Theorem 2.4 we can analyze in a more
systematic way the behavior of $d$-span of finite (and general)
sets $Z$ in terms of the mutual distances between the points of
$Z$. This can be done in terms of the so-called $\beta$-spread,
introduced in \cite{Yom1,Yom2}. A very closely related notion is
the ``$\beta$-weight of minimal spanning trees" (see
\cite{Koz.Lot.Stu} and references therein). Some results in this
direction are given in Section 3.2.6 below.

\subsection{Examples in one dimension and the Favard bound}

We complete the present section with writing down explicitly the
resulting bounds for some one-dimensional sets $Z$. Let us start
with a regular grid. \bt Let $d$ be given and let
$G^1_s=\{x_1=-1,x_2,\dots,x_s=1\}$ be a regular grid in $[-1,1], \
s>d$. Then $R_d(G^1_s)\leq T_d({{4-\mu}\over {\mu}})$ where
$\mu={{2(s-d)}\over {(s-1)}}$. In particular, $R_d(G^1_s)$ is
finite for $s>d$ and it tends to $1$ for $s\rightarrow \infty$.\et
\pr It follows from the bounds on $\o_d(G^1_s)$ computed in
Example 1, Section 3. The result of Theorem 2.5 was obtained by a
different method in \cite{Zah,Yom.Zah}. $\blacksquare$

\medskip

In Section 3.1 below we compute $\o_d(Z)$ for $Z=Z_r=\{1,{1\over
{2^r}},{1\over {3^r}},\dots,{1\over {k^r}},\dots\}$. We get
$$\o_d(Z_r)\asymp {{r^r}\over {(r+1)^{r+1}}} \
{1\over {d^r}}.$$ In particular, for $r=1$ i.e for
$Z_1=\{1,{1\over 2},{1\over 3},\dots,{1\over k},\dots\}$ we get
$\o_d(Z_1) \asymp {1\over {4d}}.$

\medskip

Now for for $0< q < 1$ and for
$Z(q)=\{1,q,q^2,q^3,\dots,q^m,\dots\}$ computations in Section 3.1
give $\o_d(Z(q)) \asymp {{q^d}\over {\log ({1\over q})}}$.

Substituting these expressions for the $d$-span into the
expression of Theorem 1.3 we obtain: \bt For the sets $Z_r$ and
$Z(q)$ as above $$R_d(Z_r)\leq R_d(\o)= T_d ({{2-\o}\over \o}),$$
where $\o \asymp {{r^r}\over {(r+1)^{r+1}}} \ {1\over {d^r}}$ or
$\o \asymp {{q^d}\over {\log ({1\over q})}}$, respectively. For
$d\rightarrow \infty$ we have $R_d(\o)\asymp ({{4(r+1)^{r+1}}\over
{r^r}})^d d^{rd}$ or $R_d(\o)\asymp (4\log {1\over q})^d ({1\over
q})^{d^2},$ respectively. In particular, $$R_d(Z_1)\preceq
2^{4d}d^d, \ R_d(Z({1\over 2}))\preceq 2^{d^2+2d}.$$\et In
\cite{Fav} slightly better bounds are given in the last two
examples:
$$R_d(Z_1)\leq (2d)^d, \ R_d(Z({1\over 2}))\leq
(d+1)2^{{d^2+3d-2}\over 2}.$$ Favard's method for bounding
$R_d(Z)$ is to fix $d+1$ points $Z$ and to estimate the
corresponding interpolation polynomial of degree $d$. This
produces the following general bound (\cite{Fav}): \be R_d(Z)\leq
\inf_{x_1,\dots,x_{d+1}\in Z}\sum_{i=1}^{d+1}{1\over
{A'_{x_1,\dots,x_{d+1}}(x_i)}},\ee where
$A_{x_1,\dots,x_{d+1}}(x)=\prod_{j=1}^{d+1}(x-x_j).$

\smallskip

Unfortunately, we cannot expect Favard's approach to
produce realistic bounds on $R_d(Z)$ for general one-dimensional $Z\subset
[-1,1]$. The problem is that considering polynomials of degree $d$
we analyze the finite subsets in $Z$ containing exactly $d+1$
points, and therefore we cannot take into account the influence of
the rest of the set $Z$. In the examples $Z_r$ and $Z(q)$
considered above this method works well since for each $d$ the
first $d+1$ points of these sets give a sufficiently accurate
approximation of the entire set. However, for a uniform grid
$G^1_s, \ s\gg d,$ a straightforward application of the Favard
estimate gives $R_d(G^1_s)\leq (2e)^d$ (the minimum in (2.2) being
achieved on the approximately uniform sub-grid formed by $d+1$
points in $G^1_s$), and this bound does not depend on $s$ at all.
Our bound given by Theorem 2.6 (which in this case is sharp up to
a constant) shows that $R_d(G^1_s)$ for any fixed degree $d$
indeed tends to $1$ as $s$ increases.

\medskip

\noindent {\bf Remark.} It is an interesting problem to
investigate the asymptotic behavior of $R_d(Z)$ as $d\rightarrow
\infty$ for ``fractal" sets $Z$ in one and several dimensions. The
examples above give some hope that the metric $d$-spread, being a
rather coarse metric invariant, still provides an adequate tool
for this problem. On the other hand, as it was mentioned above,
there are $d$ -definite sets $Z$ for which $\o_d(Z)=0$.

As for a regular grid $G^n_s$ with the step ${1\over s}$ in the
unit cube $Q^n_1$ we notice that the following inequality is true:
\bl For each n $R_d(G^n_s)\leq [R_d(G^1_s)]^n$.\el \pr Induction
by the dimension.

\section{More examples of $d$-definite sets} \setcounter{equation}{0}

In this section we consider in somewhat more details properties of the
$d$-span and present more examples, stressing the question of positivity of
$\omega_d(Z)$. In particular, we provide the proofs of Theorems 2.2 and 2.3
and of some results used in Section 2.2 above.

\subsection{Some one-dimensional examples}

For $n=1$ the sub-level set $V=V_\rho(P)\subset [-1,1]$ is just a
finite union of closed intervals. The maximal possible number of
these intervals is $d=deg P$. Clearly, the covering number $M(\e,V)$
satisfies $M(\e,V)\leq d+\mu_1(V){1\over \e},$ in agreement with
Theorem 1.2 above. We get \bp For a set $Z\subset [-1,1]$,
$\o_d(Z)=\sup_{\e>0} \e(M(\e,Z)-d).$\ep This immediately implies
\bc For $\vert Z \vert \leq d$ we have $\o_d(Z) = 0$. For
$\vert Z \vert > d$ the $d$-span $\o_d(Z)$ is strictly
positive.\ec In fact, the following more accurate bound can be
given: \bp Let $\vert Z \vert = p > d$ and let $\e_0$ be the
minimal distance between the points of $Z$. Then the $d$-span
$\o_d(Z)$ satisfies the inequality $\o_d(Z)\geq \e_0(p-d).$\ep
\pr We have $M(\e_0,Z)=p$. $\blacksquare$

\smallskip

In Section 3.3.5 below we generalize this last remark to higher
dimensions.

\medskip

Let us give now some initial specific examples where the $d$-span
can be explicitly estimated.

\medskip

\noindent{\bf Example 1.} Let $G_s=\{x_1=-1,x_2,\dots,x_s=1\}$ be
a regular grid in $[-1,1]$. The covering number $M(\e,G_s)$ is
$[{2\over \e}]+1$ for $\e \geq {2\over {s-1}}$, and it is $s$ for
$\e < {2\over {s-1}}$. Therefore the function $\e(M(\e,G_s)-d)$
behaves as $2-d\e$ for $\e \geq {2\over {s-1}}$, and it is
$\e(s-d)$ for $\e < {2\over {s-1}}$. As Corollary 3 above shows,
for $s \leq d$ we get $\o_d(G_s)= 0$. For $s>d$ the supremum is
achieved for $\e={2\over {s-1}}$ and we get $\o_d(G_s)=
{{2(s-d)}\over {s-1}}$. Notice that $\o_d(G_s)$ tends to the total
length of $[-1,1]$ as $s$ grows (or as the ``density" of the set
$G_s$ inside $[-1,1]$ increases).

\medskip

\noindent{\bf Example 2.} Let $Z_r=\{1,{1\over {2^r}},{1\over
{3^r}},\dots,{1\over {k^r}},\dots\}$. An easy computation shows
that $M(\e,Z_r)\asymp ({1\over \e})^{1\over {r+1}}.$ Hence
$$\o_d(Z_r)=\sup_{\e>0} \e(M(\e,Z_r)-d)\asymp \sup_{\e>0} \e(({1\over
\e})^{1\over {r+1}}-d)={{r^r}\over {(r+1)^{r+1}}} \ {1\over
{d^r}},$$ the supremum being attained for $\e=({r\over {r+1}})^r \
{1\over {d^{r+1}}}$. In particular, for $Z_1=\{1,{1\over
2},{1\over 3},\dots,{1\over k},\dots\}$ we get $\o_d(Z_1) \asymp
{1\over {4d}},$ the supremum being attained for $\e={1\over
{4d^2}}$.

\medskip

\noindent{\bf Example 3.} Let for $0< q < 1,$ \
$Z(q)=\{1,q,q^2,q^3,\dots,q^m,\dots\}$. Computations as above give
$\o_d(Z(q)) \asymp {{q^d}\over {\log ({1\over q})}}$.

\medskip

As for sharpness of these bounds, we show it via Theorem 2.1
above, which claims that for a sublevel set $V=V_\rho(P)$
containing $Z$ we have $\mu(V)\geq \o_D(Z)$. Now the sets of the
form $V_\rho(P)$ are exactly all the sets containing at most $d$
intervals. Therefore if we can cover $Z$ by $d$ intervals of a
total length $a$ then by Theorem 2.1 we have $\o_D(Z)\leq a$.

\smallskip

In Example 1 let us cover the grid $G_s$ with $d$ intervals, each
containing $[{s\over d}]$ consecutive points. There are $d-1$ gaps
of the length ${2\over {s-1}}$ between these intervals, so their
total length is $2-{{2(d-1)}\over {s-1}}={{2(s-d)}\over {s-1}}$.
So the bound above is sharp.

\smallskip

In Example 2 we can easily find a covering of the set $Z_r$ with
$d$ intervals of the total length ${1\over {d^r}}$. Indeed, take
first $d-1$ intervals of a small length, each covering exactly one
point from $\{1,{1\over {2^r}},{1\over {3^r}},\dots,{1\over
{({d-1})^r}}\}.$ The rest of the set $Z_r$ we cover by one
interval of the length ${1\over {d^r}}$. So also here the bound
above is sharp, up to a constant.

The same is true also in Example 3.

\subsection{Higher dimensions}

Let us start with some simple general properties of the $d$-span.
Certainly, this geometric invariant is ``stronger" than the usual
$n$-measure $\mu_n$: \bp For a measurable subset $Z\subset B^n_1$
the $d$-span $\omega_d(Z)$ satisfies $\omega_d(Z)\geq
\mu_n(Z).$\ep \pr Take $\e \rightarrow 0$ in Definition 1.2,
notice that $M_d(\e)$ grows at most as $({1\over \e})^{n-1}$, and
use the fact that if we can cover $Z$ by $M(\e,Z)$ disjoint
$\e$-cubes then $M(\e,Z)\geq {{\mu(Z)}\over {\e}^n}. \ \
\blacksquare$

\subsubsection{Sets of positive $s$-Hausdorff measure, $n-1<s<n$}

The result above can be generalized to sets of fractal Hausdorff
measures. Let us recall that for $\b>0$ the $\b$-Hausdorff measure
of $Z$ is defined as $$ H_\b(Z) = \lim_{\a \rightarrow
0}H^\a_\b(Z), $$ where $H^\a_\b(Z)$ is the lower bound of all the
sums of the form $\sum^\infty_{i=1} r^\b_i, \ r_i\ \leq\ \a$ and
$Z\subset \cup^\infty_{i=1}A_i$, with the $diam\ A_i\ \leq\ r_i$.
(See e.g. \cite{Fal}).

\medskip

However, in case $s<n$ we need more geometric information on our
set $Z$ (and not only the positivity of its $s$-Hausdorff measure
$H_s(Z)$) to conclude that the volume of any simple semi-algebraic
set containing $Z$ is large. Indeed, think about a long but
rapidly oscillating curve inside a small ball in the plane.

What we need is a kind of an ``injectivity radius" $\e^0$ of $Z$
for which the covering $\e^0$ balls are almost disjoint. Let us
give the following definition: \bd Let $H_s(Z)>0$. We define the
$s$-injectivity radius $\a^0_s(Z)$ as the maximal $\a$ such that
$H^{\a'}_s(Z)\geq {1\over 2}H_s(Z)$ for all $\a'\leq \a.$ \ed Now
we can compare the covering number and the $s$-Hausdorff measure:
\bp For $\e\leq \hat \a={1\over {\sqrt n}}\a^0_s(Z)$ we have
$M(\e,Z)\geq {1\over {2\sqrt {n^s}}}H_s(Z)({1\over \e})^s.$\ep \pr
By definition of $\a^0_s(Z)$ and of $H^\a_s(Z)$ we have for any
covering of $Z$ by $M(\e,Z)$ $\e$-cubes $$M(\e,Z)(\sqrt n \e)^s
\geq {1\over 2}H_s(Z).$$ Hence $M(\e,Z)\geq {1\over {2\sqrt
{n^s}}}H_s(Z)({1\over \e})^s.\ \ \blacksquare$

\medskip

Let us fix $n$ and $d$. We fix also a certain $s=n-1+\sigma, \
\sigma>0$. We can prove now a general lower bound for the $d$-span
of sets $Z$ with positive Hausdorff $s$-measure.

\smallskip

Let us introduce some notations. As above, we have
$M_d(\e)=C_0(n,d)+ C_1(n,d)({1\over{\e}}) + \dots +
C_{n-1}(n,d)({1\over{\e}})^{n-1},$ where the constants
$C_0(n,d),\dots,C_{n-1}(n,d)$ depending only on $n,d$ have been defined
in Theorem 1.2 above. For small $\e$ the leading term of
degree $n-1$ in ${1\over \e}$ in $M_d(\e)$ determines the
asymptotic behavior of this expression, so let us define
$\e_1=\e_1(d)$ as the maximal $\e$ such that $M_d(\e')\leq
2C_{n-1}(n,d)({1\over{\e'}})^{n-1}$ for all $\e'\leq \e$. Finally, for any $H>0$
let us put $\e_2=\e_2(H,d)=[{H\over {8C_{n-1}(n,d)\sqrt
{n^s}}}]^{1\over \sigma}.$ Now we are ready to state the result.
\bt Let $s=n-1+\sigma, \ \sigma>0$, and let $Z\subset Q^n_1$
satisfy $H_s(Z)=H>0$. Then $$\o_d(Z)\geq {1\over 4}{\hat
\e}^{1-\sigma}H_s(Z).$$ Here $\hat \e = min \{\hat
\a,\e_1(d),\e_2(H,d)\}$.\et \pr By definition \be
\o_d(Z)=\sup_{\e>0} \ e[M(Z,\e)-M_d(\e)]\geq \hat e[M(Z,\hat
\e)-M_d(\hat \e)].\ee By the choice of $\hat \e$ and by
Proposition 3.4 we have $M(\hat \e,Z)\geq {1\over {2\sqrt
{n^s}}}H_s(Z)({1\over \hat \e})^s,$ while $M_d(\hat \e)\leq
2C_{n-1}(n,d)({1\over{\hat \e}})^{n-1}$ since $\hat \e \leq
\e_1(d).$ Therefore \be M(Z,\hat \e)-M_D(\hat \e)\geq {1\over
{2\sqrt {n^s}}}H_s(Z)({1\over \hat \e})^s -
2C_{n-1}(n,d)({1\over{\hat \e}})^{n-1}.\ee Finally, the condition
that $\hat \e \leq \e_2(H,d)=[{H\over {8C_{n-1}(n,d)\sqrt
{n^s}}}]^{1\over \sigma}$ implies that the right-hand side of
(3.2) is not smaller than ${1\over {4\sqrt {n^s}}}H_s(Z)({1\over
\hat \e})^s$. Combining this last inequality with (3.1) and (3.2)
we obtain the required bound. $\blacksquare$

\medskip

\noindent {\bf Remark.} An important feature of Theorem 3.1 is
that we do not need to assume that the $s$-Hausdorff measure of
$Z$ is ``large". Just the fact that $H_s(Z)>0$ implies
$\o_d(Z)>0$. To stress the dependence of the bound of Theorem 3.1
on $s$ and $H_s(Z)$ let us assume that the radius of injectivity
$\a^0_s(Z)$ is large while the measure $H=H_s(Z)$ is small. Then
$\hat \e = min \{\hat \a, \e_1(d),\e_2(H,d)\}=\e_2(H,d)=\tilde
C(d,s)H^{1\over \sigma},$ and therefore by Theorem 3.1 \
$\o_d(Z)\sim [C_1H]^{1\over \sigma}.$ This bound blows up as
$\sigma \rightarrow 0$ or $s \rightarrow n-1.$

Being quite effective, the bound of Theorem 3.1 is not sharp.
Compare \cite{Bru,Bru1}.

\smallskip

However, for $s$ exactly equal to $n-1$ there is still a
possibility to bound $\o_d(Z)$ from below if $H_s(Z)$ is strictly
greater than $2\sqrt{n^{n-1}}C_{n-1}(n,d)$. This bound is obtained
in Corollary 4 in Section 3.2.4 below.

\subsubsection{Sets with large covering number}

The following result is parallel to Theorem 3.1, but it replaces
the assumption of positivity of $H_s(Z)$ with the assumption that
the covering number $M(Z,\e)$ grows as $C_s({1\over \e})^s, \
s>n-1$, for $\e$ sufficiently small. We preserve essentially the
same notation as in Theorem 3.1: define the s-covering injectivity
radius $\e^0_s(Z)$ as the maximal $\e$ such that $M(Z,\e')\geq
{1\over 2}C_s({1\over \e'})^s$ for all $\e'\leq \e$. The parameter
$\e_1=\e_1(d)$ is defined exactly as above, and we put
$e'_2(C_s,d)=[{C_s\over {C_{n-1}(n,d)}}]^{1\over s}$. \bt Let
$s=n-1+\sigma, \ \sigma>0$, and let $Z\subset Q^n_1$ satisfy
$M(Z,\e)\geq C_s({1\over \e})^s$, for all sufficiently small $\e$.
Then $$\o_d(Z)\geq {1\over 4}{\hat \e}^{1-\sigma}C_s.$$ Here $\hat
\e = min \{\e^0_s(Z),\e_1(d),\e'_2(C_s,d)\}$.\et \pr Exactly the
same as for Theorem 3.1. $\blacksquare$

\medskip

\noindent {\bf Remark.} As above, if $\hat \e = \
e'_2(C_s,d)=[{C_s\over {C_{n-1}(n,d)}}]^{1\over s}$ then we get
$$\o_d(Z)\geq {1\over 4}({1\over {8C_{n-1}(n,d)}})^{{1-\sigma}\over
\sigma}C_s^{1\over \sigma}.$$ Notice also that Theorem 3.2
formally implies Theorem 3.1 because of Proposition 3.4. However,
since the Hausdorff measure is probably a somewhat more natural
invariant than the covering number, it looks preferable to
separate these two statements.

\subsubsection{Entropy and Hausdorff dimension}

We recall here the notions of the entropy and the Hausdorff
dimensions.

\bd Let $A\subset X$ be a bounded subset in a metric space $X$.
\item{1.} $\dim_eA=\inf\{\b, \ \exists K,$ such that for each $\e>0, \
\ N(\e,A)\ \leq\ K({1\over\e})^\b\}$ is called the entropy
dimension of $A$.

\smallskip

\item{2.} $\dim_H A=\inf\{\b, \ S_\b(A)<\infty\}$ is called the
Hausdorff dimension of $A$.\ed The notion of the entropy dimension
appears in fractal geometry under many different names, in
particular: ``Minkowski dimension" - probably, the most justified
historically, - ``capacity dimension", ``box dimension".

\medskip

It is well known (see, for example, \cite{Fal}) that for any set
$A$ we have $\dim_e A\geq \dim_H A$. In particular, for countable
sets $A$ always $\dim_H A=0$ while $\dim_e A$ may take any value.
The bounds of Theorems 3.1 and 3.2 imply the following: \bp For
any $d$ and for any
subset $Z\subset B^n_1$ if $\dim_e Z > n-1$ then $\o_d(Z)>0$. In
particular, this is true if \ \ \ $\dim_H Z > n-1$.\ep

\subsubsection{Sets of dimension $n-1$}

Now we consider the case $s=n-1$. Here we start with the covering
number and obtain the corresponding result for the Hausdorff
measure as a corollary. Let $M(Z,\e)\geq C({1\over \e})^{n-1}$ for
$\e \leq \e^0_{n-1}$, with $C>C_{n-1}(n,d)$. We define $\e'_1(d,C)$
as the largest $\e$ for which $M_d(\e'')\leq Q(C,d)({1\over
\e''})^{n-1}$ for all $\e''\leq \e$. Here
$Q(C,d)=[C_{n-1}(n,d)+{1\over 2}(C-C_{n-1}(n,d))].$ \bt Let $Z\subset
Q^n_1$ satisfy $M(Z,\e)\geq C({1\over \e})^{n-1}$ for $\e \leq
\e^0_{n-1}$, with $C>C_{n-1}(n,d)$. Then
$$\o_d(Z)\geq {1\over 2}{\hat \e}(C-C_{n-1}(n,d)),$$ where $\hat \e
= min \{\e^0_{n-1},\e'_1(C,d)\}$.\et \pr Exactly as in Theorem
3.1. $\blacksquare$

\medskip

Via Proposition 3.4 we obtain: \bc Let $Z\subset B^n_1$ satisfy
$H_{n-1}(Z)>2\sqrt{n^{n-1}}C_{n-1}(n,d).$ Then
$$\o_d(Z)\geq {1\over 2}{\hat \e}(C-C_{n-1}(n,d)),$$ where
$C={1\over {2\sqrt{n^{n-1}}}}H_{n-1}(Z)$ and $\hat \e = \min \{\
{1\over {\sqrt n}}\alpha^0_{n-1}(Z),\e'_1(C,d)\}$.\ec

\subsubsection{Dense finite subsets in ``massive" sets}

Each of the results above produces, in particular, a finite subset
$Z'\subset Q^n_1$ with $\o_d(Z')>0$. Indeed, in each of the
situations covered by Theorems 3.1-3.3 and Corollary 4 let us
define $Z'\subset Z$ as the set of the centers of all the ${1\over
2}\hat \e$-cubes providing a covering of $Z$ with $M(Z,{1\over
2}\hat \e)$ elements. We have \bt In each of the situations
covered by Theorems 3.1-3.3 and Corollary 4 the $d$-span of the
finite set $Z'\subset Z$ satisfies $\o_d(Z')\geq ({1\over 2})^n K$
where $K$ is the appropriate lower bound for $\o_d(Z)$.\et \pr If
certain ${1\over 2}\hat \e$-cubes cover $Z'$ then the
corresponding $\hat \e$-cubes cover $Z$. Therefore $M(Z',{1\over
2}\hat \e)\geq M(Z,{1\over 2}\hat \e)$. The rest of the proof goes
exactly as in the results above. $\blacksquare$

\subsubsection{Bounding $D$-span via Minimal spanning trees}

Theorem 3.4 provides a class of examples of finite subsets
$Z\subset B^n_1$ with positive $d$-span: roughly, those are
sufficiently dense finite subsets of sets of dimension $n-1$ or
higher. It is important to analyze the behavior of $d$-span of
finite (and general) sets $Z$, given by themselves, with no
relation to an underlying ``large" set, in terms of their metric
structure. Here an appropriate invariant may be the so-called
$\beta$-spread, introduced in \cite{Yom1,Yom2}. A very closely
related notion is the ``$\beta$-weight of minimal spanning trees"
(see \cite{Koz.Lot.Stu} and references therein). The main reason
for us to relate the $d$-span with the $\beta$-spread and minimal
spanning trees is that a lot of information is available today in
this direction (for some initial references see
\cite{Koz.Lot.Stu}), and we can hope to ultimately incorporate
this information in our study of polynomial and smooth
interpolation problems.

\smallskip

Let's recall a definition of $\beta$-spread. Let $G_p$ be the set
of all connected non-oriented trees with $p$ vertices. We write
$(i,j)\in g$, for $g\in G_p$, if the vertices $i$ and $j$ are
connected by the edge in $g$.

\bd Let $X$ be a metric space, $\b>0$. For each $x_1,
\cdots,x_p\in X$ and $g\in G_p$ let
$\rho_\b(g,x_1,\cdots,x_p)=\sum_ {(i,j)\in g} d(x_i,x_j)^\b$,
where $d$ is a distance in $X$. Define $\rho_\b (x_1,\cdots,x_p)$
as $\inf_{g\in G_p} \rho_\b(g,x_1,\cdots,x_p)$. The tree $g$ on
which the infinum is achieved is called the $\b$-minimal spanning
tree. Now let $A\subset X$. We define the $\b$-spread of $A, \
V_\b(A)$, by
$$ V_\b(A)=\sup_{p,x_1,\cdots,x_p\in A}\rho_\b(x_1,\cdots,x_p).
$$\ed For $x_1,\cdots,x_p\in X$, \ $\rho_\b (x_1,\cdots,x_p)$ is called a
$\beta$-weight of the minimal spanning tree $g$ on
$(x_1,\cdots,x_p)$. Notice that the $1$-minimal tree is also
minimal for any $\b$ (see \cite{Koz.Lot.Stu}).

\smallskip

Under a different name $\b$-spread for subsets of a real line has
been studied in \cite{Bes.Tay}. A notion of $\beta$-weight has
appeared earlier in geometric combinatorics and in fractal
geometry. Compare \cite{Koz.Lot.Stu,Lap.Fra}, \cite{Fal} and
references therein. However, we are not aware of any appearance of
the general notion of $\b$-spread in metric spaces, as defined
above.

Let us also notice that as a function of $\b$ the spread $V_\b(A)$
is a kind of a zeta-function. For $A=\{0,1,3,6,10,15,\dots,
{1\over 2}n(n+1),\dots \}$ the spread $V_\b(A)$ is exactly the
Riemann $\zeta$-function $\zeta(-\beta)$, while for
$A=\{\alpha_0,\alpha_1,\dots,\alpha_n,\dots\}$ with $\alpha_0=0, \
\alpha_1=1,\dots,\alpha_n=\sum^n_{i=1}{1\over i}, \ n\geq 1,$  we
have $V_\b(A)=\zeta(\beta)$. So it may be a good idea to
substitute into $V_\b(A)$ complex values of $\b$. See
\cite{Lap.Fra} for a detailed treatment of fractal geometry from
this point of view.

\medskip

We shall not touch here general properties of $\beta$-spread, as
well as its relations to the geometry of critical values of smooth
functions. Instead we give a lower bound for the $d$-span in terms
of $\beta$-spread. Let us first provide an immediate
generalization to higher dimensions of Proposition 3.2 above. We
have to consider here the $l^\infty$ distance instead of the usual
Euclidean distance in ${\mathbb R}^n$. \bp Let $d$ and a
finite subset $Z\subset Q^n_1$ be given, and let $\e_0$ be the
minimal distance between the points of $Z$. Assume that $\vert Z
\vert = p
> M_d(\e_0).$ Then the $d$-span $\o_d(Z)$ satisfies the inequality
$\o_d(Z)\geq \e^n_0(p-M_d(\e_0))>0.$\ep For any $Z$, finite or
infinite, we can apply Proposition 3.6 to finite subsets of $Z$.
Let us introduce some convenient notations (see \cite{Yom2}). \bd
Let $X$ be a metric space. For $x_1,\cdots,x_p\in X$, let
$\nu(x_1,\cdots,x_p) = \min_{i\neq j} \ d(x_i,x_j)$. For $A\subset
X$ define $\eta_A(p)$ for any natural $p\ \geq\ 2$ by
$$ \eta_A(p)=\sup_{x_1,\cdots,x_p\in A}\ \nu(x_1,\cdots,x_p).$$\ed
Proposition 3.6 implies \bc Let $d$ and a finite subset $Z\subset
Q^n_1$ be given. If for a certain $p$ we have $M_d(\eta_Z(p))<p$
then $\o_d(Z)\geq \eta^n_Z(p)(p-M_d(\eta_Z(p)))>0$.\ec Let us
remind that $M_d(\e)=\sum_{j=0}^{n-1}C_j(n,d)({1\over \e})^j \leq
C'(n,d)({1\over \e})^{n-1}$ for $\e\leq 1$. Hence we have a weaker
but simpler version of Corollary 5: \bc Let $d$ and a finite
subset $Z\subset Q^n_1$ be given. If for a certain $p$ we have
$C'(n,d)(\eta_Z(p))^{1-n}<p$ then \be \o_d(Z)\geq
\eta^n_Z(p)(p-C'(n,d)(\eta_Z(p))^{1-n})>0.\ee \ec Now we can give
a criterion of positivity of $\o_d(Z)$ in terms of the
$\beta$-spread of $Z$: \bt Let $d$ and a subset $Z\subset Q^n_1$
be given. If for a certain $\beta,$ with $n-1<\beta \leq n,$ we
have $V_\beta(Z)>C'(d)^{\beta\over {n-1}}\zeta({\beta\over
{n-1}}),$ where $\zeta(x)=\sum_{p=1}^\infty {1\over {p^x}},$ then
$\o_d(Z)>0$.\et \pr Assume that $\o_d(Z)=0$. Then by Corollary 5
we have for each $p=2,3,\dots$ that $C'(n,d)(\eta_Z(p))^{1-n}\geq
p$. Hence $\eta_Z(p)\leq ({{C'(n,d)}\over p})^{1\over {n-1}}$ and
for each $\beta$ we have \be \sum_{p=2}^\infty \eta^\beta_Z(p)
\leq (C'(n,d))^{\beta\over {n-1}}\sum_{p=1}^\infty ({1\over
p})^{\beta\over {n-1}}=C'(n,d)^{\beta\over {n-1}}\zeta({\beta\over
{n-1}}).\ee Now, the following result relates $\eta_Z(p)$ and
$V_\b(Z)$: \bp For any $\b>0$
$$ \sup_{p\geq 2}(p-1)\eta^\b_Z(p)\ \leq\ V_\b(Z)\ \leq\ \sum^\infty_{j=2}
\eta^\b_Z(j). $$ \ep The proof of Proposition 3.7 is given in
\cite{Yom1} (see also \cite{Yom2}). Combining this result with
(3.4) we complete the proof of Theorem 3.5. $\blacksquare$

\medskip

In analogy with the Hausdorf and entropy dimensions let us define
the $V$-dimension as follows: $\dim_V A =\inf\{\b, \
V_\b(A)<\infty\}$. It turns out that always $\dim_V A=\dim_e A$
(see \cite{Yom1,Yom2,Koz.Lot.Stu}). Now Theorem 3.5 provides
another proof of Proposition 3.5 above. Indeed, if $\dim_e Z=s
>n-1,$ we fix some $\beta$ such that $n-1 < \beta < s = \dim_V Z$. By
definition of $\dim_V Z$ we have $V_\beta(Z)=\infty$ while
$\zeta({\beta\over {n-1}})$ is finite since $\beta>n-1$. Theorem
3.5 implies now that $\o_d(Z)>0$.

\smallskip

There are limit cases where $\beta$-spread is more sensitive to
certain subtle geometric properties of $Z$ than the covering
number (see \cite{Yom1,Yom2,Koz.Lot.Stu} and references therein).
It is also related with some important notions in Potential
Theory, like transfinite diameter. We plan to present some results
in this direction separately.

\bibliographystyle{amsplain}

\begin{thebibliography}{10}

\bibitem{Ber1} S. Bernstein, Sur une formule d'interpolation, {\sl
Comptes Rendus AS} {\bf 191} (1930), 635-637.

\bibitem{Ber2} S. Bernstein, Sur la limitation des valeurs d'un
polynome $P_n(x)$ de degr\'e $n$ sur tout un segment par ses
vleurs en $n+1$ points du segment, {\sl Isvestiya AN SSSR},
(1931), 1025-1050.

\bibitem{Ber3} S. Bernstein, On the trigonometric interpolation
via the least square method, {\sl Doc. AN SSSR}, {\bf 4} (1934),
1-8.

\bibitem{Bes.Tay} A. S. Besicovitch, J. Taylor, On the
complementary intervals of a linear closed set of zero Lebesgue
measure, {\sl J. London Math. Soc.} {\bf 29} (1954), 449-459.

\bibitem{Bie.Mil.Paw} E. Bierstone; P. Milman, W. Pawlucki,
Differentiable functions defined in closed sets. A problem of
Whitney, {\sl Invent. Math.} {\bf 151} (2003), no. 2, 329-352.

\bibitem{Bor.Naz.Sod} A. Borichev, F. Nazarov, M. Sodin, Lower
bounds for quasianalytic functions. II. The Bernstein
quasianalytic functions. {\sl Math. Scand.} {\bf 95} (2004), no.
1, 44--58.

\bibitem{Bru} A. Brudnyi, On Covering Numbers of Sublevels Sets of Analytic
Functions, to appear in J. of Appr. Theory.

\bibitem{Bru1} A. Brudnyi, On a BMO-property for subharmonic
functions. J. Fourier Anal. Appl. 8 (2002), no. 6, 603--612.

\bibitem{Bru.Bru2} A. Brudnyi, Yu. Brudnyi, Remez Type Inequalities
and Morrey-Campanato Spaces on Ahlfors Regular Sets, Contemporary
Mathematics, {\bf 445} (2007), 19-44.

\bibitem{Bru.Gan} Yu. Brudnyi, M. Ganzburg, On an extremal problem
for polynomials of $n$ variables, {\sl Math. USSR Izv.} {\bf 37}
(1973), 344-355.

\bibitem{Bru.Shv} Y. Brudnyi, P. Shvartsman, Whitney's extension
problem for multivariate $C^{1,\o}$-functions, {\sl Trans. Amer.
Math. Soc.} {\bf 353} (2001), 2487-2512.

\bibitem{Cop.Riv} D. Coppersmith, T. J. Rivlin, The growth of
polynomials bounded at equally spaced points, {\sl SIAM J. Math.
Anal.} {\bf 23} (1992), no. 4, 970--983.

\bibitem{Erd} T. Erdelyi, Remez-type inequalities and their
applications, {\sl J. Comp. Appl. Math.} {\bf 47} (1993) 167-209.

\bibitem{Fav} J. Favard, Sur l'interpolation, {\sl Bull. de la S.
M. F.}, {\bf 67} (1939), 103-113.

\bibitem{Fef} Ch. Fefferman, Whitney's extension problem for $C\sp
m$, {\sl Ann. of Math.} (2) {\bf 164} (2006), no. 1, 313--359.

\bibitem{For} T. Fort, Finite differences and difference relations
in the real domain, {\it Oxford, at the Clarendon Press}, 1948.

\bibitem{Gor} E. A. Gorin, A. Cartan's lemma following B. Ya.
Levin with applications, {\sl J. of Math. Phys., Anal., Geom.,}
{\bf 3}, 1 (2007), 13-38 (Russian).

\bibitem{Hil} F. B. Hildebrand, Introduction to numerical analysis,
Second Edition, {\it Dover Publications, Inc}, New-York, 1987

\bibitem{Fal} K. Falconer, Fractal geometry. Mathematical foundations
and applications, Second edition. John Wiley and Sons, Inc.,
Hoboken, NJ, 2003. xxviii+337 pp.

\bibitem{Iva}  L. D. Ivanov,
Variazii mno\^zhestv i funktsii. (Russian) [Variations of sets and
functions] Edited by A. G. Vituskin, {\sl Izdat. ``Nauka''},
Moscow, (1975), 352~p.

\bibitem{Kol} A. N. Kolmogorov, Asymptotic characteristics of some
completely bounded metric spaces, {\sl Dokl. Akad. Nauk SSSR} {\bf
108} (1956), 585-589.

\bibitem{Kol.Tih} A. N. Kolmogorov; V. M. Tihomirov,
$\varepsilon $-entropy and $\varepsilon $-capacity of sets in
functional space. {\sl Amer. Math. Soc. Transl.} {\bf  17},
(1961), 277-364.

\bibitem{Koz.Lot.Stu} G. Kozma, Z. Lotker, G. Stupp, The
minimal spanning tree and the upper box dimension, {\sl Proc.
Amer. Math. Soc.} {\bf 134} (2006), no. 4, 1183--1187.

\bibitem{Lap.Fra} M. Lapidus, M. van Frankenhuysen, Fractal
Geometry and Number Theory, Birkhauser, 2000.

\bibitem{Min} H. Minkowski, Theorie der konvexen Korper,
insbesondere Begrundung ihres Oberflachenbegriffs, in: Gesammelte
Abhandlungen von Hermann Minkowski (part II, Chapter XXV),
Chelsea, New York, 1967, 131-229.

\bibitem{Nad} B. Nadler, private communication.

\bibitem{Lor} R. Lorentz, Multivariate Birkhoff Interpolation,
{\sl Schpringer-Verlag}, 1992.

\bibitem{Naz.Sod.Vol1} F. Nazarov, M. Sodin, A. Volberg,
The geometric Kannan-Lovász-Simonovits lemma, dimension-free
estimates for the distribution of the values of polynomials, and
the distribution of the zeros of random analytic functions, {\sl
Algebra i Analiz} {\bf 14} (2002), no. 2, 214--234; translation in
{\sl St. Petersburg Math. J.} {\bf 14} (2003), no. 2, 351--366.

\bibitem{Naz.Sod.Vol2} F. Nazarov, M. Sodin, A. Volberg,
Local dimension-free estimates for volumes of sublevel sets of
analytic functions, {\sl Israel J. Math.} {\bf 133} (2003),
269--283.

\bibitem{Naz.Sod.Vol3} F. Nazarov, M. Sodin, A. Volberg, Lower bounds
for quasianalytic functions. I. How to control smooth functions,
{\sl Math. Scand.} {\bf 95} (2004), no. 1, 59--79.

\bibitem{Rah} E. A. Rakhmanov, Bounds for polynomials with a unit
discrete norm. {\sl Ann. of Math.} (2) {\bf 165} (2007), no. 1,
55--88.

\bibitem{Rem} E. J. Remez, Sur une propriete des polynomes de
Tchebycheff, {\sl Comm. Inst. Sci. Kharkov} {\bf 13} (1936) 93-95.

\bibitem{Vit1} A. G. Vitushkin, O mnogomernyh Variaziyah,
Gostehisdat, Moskow, (1955).

\bibitem{Vit2} A. G. Vitushkin, Ozenka sloznosti zadachi
tabulirovaniya, Fizmatgiz, Moskow, 1959. Translation: Theory of
the transmission and processing of information. Pergamon Press,
(1961).

\bibitem{Whi} H. Whitney, Analytic extension of differentiable
functions defined in closed sets, {\sl Trans. Amer. Math. Soc.}
{\bf 36} (1934), 63-89.

\bibitem{Yom1} Y. Yomdin, Beta-spread of sets in metric spaces
and critical values of smooth functions, Preprint, MPI Bonn,
(1983).

\bibitem{Yom2} Y. Yomdin, Beta-spread of sets in metric spaces
and critical values of smooth functions, II, to appear.

\bibitem{Yom.Zah} Y. Yomdin, G. Zahavi, High-Order discretization
of singular data, to appear.
\bibitem{Zah} G. Zahavi, Ph.D thesis, Weizmann Institute, 2007.

\bibitem{Zer1} A. Zeriahi, A minimum principle for plurisubharmonic
functions, {\sl Indiana Univ. Math. J.} {\bf 56} No. 6 (2007),
2671–2696.

\bibitem{Zer2} A. Zeriahi, Volume and capacity of sublevel sets of
a Lelong class of plurisubharmonic functions, {\sl Indiana Univ.
Math. J.} {\bf 50} No. 1 (2001), 671–703.




\end{thebibliography}

\end{document}